# Estimation of the Brownian dimension of a continuous Itô process

JEAN JACOD[1], ANTOINE LEJAY[2] and DENIS TALAY[3]

[1]*Institut de Mathématiques de Jussieu, Université Pierre et Marie Curie et CNRS (UMR 7586), 4 Place Jussieu, 75252 PARIS Cedex 05, France.*

[2]*INRIA Lorraine, IECN, Campus Scientifique, BP 239, 54506 Vandœuvre-lès-Nancy Cedex, France.*

[3]*INRIA Sophia Antipolis, 2004 Route des Lucioles, BP 93, 06902 Sophia-Antipolis, France. E-mail: Denis.Talay@sophia.inria.fr*

In this paper, we consider a $d$-dimensional continuous Itô process which is observed at $n$ regularly spaced times on a given time interval $[0, T]$. This process is driven by a multidimensional Wiener process and our aim is to provide asymptotic statistical procedures which give the minimal dimension of the driving Wiener process, which is between 0 (a pure drift) and $d$. We exhibit several different procedures, all similar to asymptotic testing hypotheses.

## 1. Introduction

In numerous applications, one chooses to model a complex dynamical phenomenon by stochastic differential equations or, more generally, by semimartingales, either because random forces excite a mechanical system, because time-dependent uncertainties disturb a deterministic trend or because one aims to reduce the dimension of a large-scale system by considering that some components contribute stochastically to the evolution of the system. Respective examples of applications are mechanical oscillators submitted to random loading, prices of financial assets and molecular dynamics.

Of course, the calibration of the model is a crucial issue. A huge literature deals with the statistics of stochastic processes, particularly of diffusion processes. Parametric and nonparametric estimators of the coefficients of stochastic differential equations have been intensively studied; see, for example, the books [6] and [7] of Prakasa Rao, in which a large number of papers are quoted and analyzed. However, somewhat astonishingly, it seems to us that most of the papers consider that the dimension of the noise is known by the observer. This hypothesis is often questionable: there is no reason to a priori fix this dimension when one observes a basket of assets or a complex mechanical structure in a random environment. Actually, the last two authors of this paper were motivated to study this question by modelling and simulation issues related to the pricing of contracts based on baskets of energy prices (see O. Bardou's thesis [3]). There was no determining financial reason to fix the Brownian motion dimension to a particular value. In addition, the interest in finding a dimension as small as possible was twofold: first, one then avoids the calibration of useless diffusion matrix components; second, practitioners need the simulation of the model, and thus the computation of contract prices and of corresponding risk measures by means of Monte Carlo simulations,to be as rapid as possible.

We thus try, in this paper, to tackle the question of estimating the Brownian dimension of an Itô process from the observation of one trajectory during a finite time interval. More precisely, we aim to build estimators which provide an 'explicative Brownian dimension' $r_B$: a model driven by an $r_B$ Brownian motion satisfyingly fits the information conveyed by the observed path, whereas increasing the Brownian dimension does not allow to fit the data any better. Stated this way, the problem is obviously ill-posed, hence our first step consists of defining a reasonable framework to develop our study.

Suppose that we observe a continuous $d$-dimensional Brownian semimartingale $X = (X^i)_{1 \leq i \leq d}$ on some space $(\Omega, \mathcal{F}, (\mathcal{F}_t), \mathbf{P})$. The observation time interval is $[0, T]$ with $T$ finite. The process $X$ is a continuous Itô process, meaning that it satisfies the following assumption.

***Hypothesis (H).*** We have

$$X_t = X_0 + \int_0^t a_s \, \mathrm{d}s + \int_0^t \sigma_s \, \mathrm{d}W_s, \tag{1}$$

where $W$ is a standard $q$-dimensional Brownian motion (BM), $a$ is a predictable $\mathbf{R}^d$-valued locally bounded process and $\sigma$ is a $d \times q$-matrix-valued adapted and cadlag process (in (1), one can replace $\sigma_s$ by the left limit $\sigma_{s-}$ so as to have a predictable integrand, if one wishes).

We say that we are in the 'pure diffusion case' when $\sigma_s = \sigma(X_s)$. We set $c_s = \sigma_s \sigma_s^\star$ (so $c_s = c(X_s)$, where $c = \sigma \sigma^\star$ in the pure diffusion case). The process $c$ takes its values in the set $\mathcal{M}_d$ of all $d \times d$ symmetric non-negative matrices. We denote by rank$(\Sigma)$ the rank of any $\Sigma \in \mathcal{M}_d$.

As is well known, the same process $X$ can be written as (1) with many different Wiener processes; namely, if $(\Pi_s)$ is a progressively measurable process taking its values in the set of $q \times q$ orthogonal matrices, then $W'_t = \int_0^t \Pi_s \, \mathrm{d}W_s$ is another $q$-dimensional Wiener process and $X$ is of the form (1) with $W'$ and $\sigma'_s = \sigma_s \Pi_s^{-1}$. The *'Brownian dimension'* $r_B$ of our model is then defined as being the smallest integer $r$ such that, after such a transformation, the last $q - r$ columns of $\sigma'_s(\omega)$ vanish (outside a $\mathbf{P}(\mathrm{d}\omega) \, \mathrm{d}s$-null set, and for $s \leq T$, of course). In this case, we can forget about the last $q - r$ components of $W'$ and, in fact, write (1) with an $r$-dimensional Wiener process. Obviously, $r_B \leq d$ always, so one could start with a model (1) with $q \leq d$ always, but it is convenient for the discussion in this paper to take $q$ arbitrary.

Our aim is to make some kind of inference on this Brownian dimension $r_B$, which is also the maximal rank of $c_s$ (up to a $\mathbf{P}(d\omega)\,ds$-null set), on the basis of the observation of the variables $X_{iT/n}$ for $i = 0, 1, \ldots, n$, where $[0, T]$ is the time interval on which the process is available. Let us make some preliminary comments, in which we refer to the "ideal" and "actual" observation schemes when one observes $X$ completely over $[0, T]$, or at times $iT/n$ only, respectively.

(1) Suppose we are in the pure diffusion case, that is, $c_s = c(X_s)$, with $c$ a continuous function, and that the range of the process is the whole of $\mathbf{R}^d$ (i.e., every open subset of $\mathbf{R}^d$ is visited by $X$ on the time interval $[0, T]$ with a positive probability). Set $r(x) := \mathrm{rank}(c(x))$ and let $A(\omega)$ be the subset of $\mathbf{R}^d$ which is visited by the path $(X(\omega)_t : t \in [0, T])$. The Brownian dimension is $r_B = \sup_{x \in \mathbf{R}^d} r(x)$, but, in the ideal scheme, we observe $R(\omega) = \sup_{x \in A(\omega)} r(x)$ and so we can only assert that $r \geq R(\omega)$. The situation is similar to what happens in the nonparametric estimation of the function $c$: in the ideal scheme, this function $c$ is known on $A(\omega)$ and hopelessly unknown on $\mathbf{R}^d \setminus A(\omega)$.

(2) More generally, the only relevant quantity we might hope to 'estimate' is the (random) maximal rank

$$R(\omega) = \sup_{s \in [0, T)} \mathrm{rank}(c_s(\omega)) \tag{2}$$

(we should take the essential supremum rather than the supremum, but the two agree since $c_s$ is right-continuous in $s$). The variable $R$ is integer-valued, so its 'estimation' is more akin to testing that $R = r$ for any particular $r \in \{0, \ldots, d\}$, although it will not be a test in the ordinary sense because $R$ is random. Note that in many models, we will have that $\mathrm{rank}(c_s(\omega))$ is independent of $s$ and $\omega$; $R$ is then non-random, but this property does not really makes the analysis any easier.

(3) In the actual scheme, we will construct an integer-valued statistics $\widehat{R}_n$ which serves as an 'estimator' for $R$. We have to somehow maximize (and evaluate) the probability that $\widehat{R}_n = R$, or perhaps this probability conditional on the value taken by $R$, or conditional on the whole path of $X$ over $[0, T]$. That is, we perform a kind of 'conditional test'.

(4) We might also take a different look at the problem. Considering the model (1), we can introduce a kind of 'distance' $\Delta_r$ between the true process $X$ and the class of all processes $X'$ of the same form, but with a diffusion coefficient $c'_s$ satisfying identically $\mathrm{rank}(c'_s) \leq r$. We then construct estimators $\Delta^n_r$ for $\Delta_r$, for all values of $r$, and decide on the basis of these $\Delta^n_r$ which Brownian dimension $r_B$ is reasonable to consider for the model. The mathematical problem is then similar to the semi-parametric estimation of a parameter in the diffusion coefficient for a discretely observed diffusion with unknown drift: here, the 'parameter' is the collection of all $\Delta_r$ and the unknown (nuisance) parameters are the processes $a_s$ and $c_s$ (or $\sigma_s$).

The paper is organized as follows. In Section 2, we explain in a more precise way the 'distances' mentioned above. Section 3 is a collection of simple linear algebra results and Section 4 contains the basic limiting results needed. In Sections 5 and 6, we then put the previous results in use to develop some statistical applications and, finally, we provide some numerical experiments in Section 7.

## 2. An instructive but non-effective approach

In this section, we measure the discrepancy between the model (1) and models of the same type but with a different Brownian dimension. We denote by $\mathcal{S}_r$ the set of all cadlag adapted $d \times q$-matrix-valued processes $\sigma'$ such that $c'_s = \sigma'_s \sigma'^\star_s$ satisfies $\mathrm{rank}(c'_s) \leq r$ a.s. for all $s$. In particular, $\mathcal{S}_0$ contains only $\sigma' \equiv 0$. With any $\sigma' \in \mathcal{S}_r$, we associate the process

$$X'_t = X_0 + \int_0^t a_s \, \mathrm{d}s + \int_0^t \sigma'_s \, \mathrm{d}W_s, \tag{3}$$

with the same $a$ and the same $W$ as in (1).

A measure of the 'distance' between the two processes $X$ and $X'$ of (1) and (3), measured on the time interval $[0, t]$, is the random variable $\Delta(X, X')_t$ defined below. In the following formula, $H$ ranges through all predictable $d$-dimensional processes with $\|H_t(\omega)\| \leq 1$ for all $(\omega, t)$, $H^\star$ is the transpose, $\langle M \rangle$ is the quadratic variation process of the semimartigale $M$ and $\bullet$ denotes stochastic integration:

$$\Delta(X, X')_t = \sup_{H \,:\, \|H\| \leq 1} \langle H^\star \bullet (X - X') \rangle_t. \tag{4}$$

This measurement of the discrepancy between $X$ and $X'$ is particularly well suited to finance, where $\mathbf{E}(\Delta(X, X')_t)$ is a measure of the difference in the $\mathcal{L}^2$ sense between the portfolio evaluations when one takes the model with $X$ or the model with $X'$. Then set

$$\Delta(r; X)_t = \inf(\Delta(X, X')_t : X' \text{is given by (3), with } \sigma' \in \mathcal{S}_r) \tag{5}$$

for the 'distance' from $X$ to the set of semimartingales with Brownian dimension not more than $r$, again on the time interval $[0, t]$.

**Remark 1.** Of course, $\Delta(X, X')_t$ is not a genuine distance, for two reasons: it does not satisfy the triangle inequality (it is rather the square of a distance) and, more important, it is random. The genuine distance (which is one of Emery's distances – see [5]) is $\sqrt{\mathbf{E}(\Delta(X, X')_t)}$, provided we identify two processes which are a.s. equal, as usual.

Also, note that the two approaches – here, where $W$ is kept fixed, and in the previous section, where $W$ may be changed into another Wiener process $W'$ – look different, but are actually the same.

The next proposition shows how to compute 'explicitly' $\Delta(r; X)_t$. We denote by $\lambda(1)_s \geq \lambda(2)_s \geq \cdots \geq \lambda(d)_s \geq 0$ the eigenvalues of the matrix $c_s$ and we set

$$L(r)_t = \int_0^t \lambda(r)_s \, \mathrm{d}s. \tag{6}$$

**Proposition 2.** *For any $r = 0, \ldots, d-1$, we have $\Delta(r; X)_t = L(r+1)_t$ and the infimum in (5) is attained.*

**Proof.** It is no restriction to suppose that $q \geq d$ (if not, we can always add independent components to $W$ and, accordingly, components to $\sigma$ which are 0). Let $J$ be the $d \times q$ matrix with $(i,i)$ entry equal to 1 when $1 \leq i \leq d$, and all other entries equal to 0. We can then find two càdlàg adapted processes $\Pi_s$ and $Q_s$, with values in the sets of $d \times d$ and $q \times q$ orthogonal matrices, respectively, and such that $\sigma_s = \Pi_s \Lambda_s^{1/2} J Q_s$, where $\Lambda_s$ is the diagonal matrix with entries $\lambda(i)_s$. Note that $c_s = \Pi_s \Lambda_s \Pi_s^\star$.

Let $I_r$ (resp., $I'_r$) be the $d \times d$ matrix with $(i,i)$ entry equal to 1 when $1 \leq i \leq r$ (resp., $r+1 \leq i \leq d$) and all other entries equal to 0. We then set $\sigma'_s = \Pi_s \Lambda_s^{1/2} I_r J Q_s$ and associate $X'$ by (3). Then $\sigma_s - \sigma'_s = \Pi_s \Lambda_s^{1/2} I'_r J Q_s$ and thus

$$\langle H^\star \bullet (X - X') \rangle_t = \int_0^t (H_s^\star \Pi_s \Lambda_s^{1/2} I'_r J Q_s Q_s^\star J^\star I'_r \Lambda_s^{1/2} \Pi_s^\star H_s) \, \mathrm{d}s.$$

The integrand above is simply $H_s^\star \Pi_s \Lambda_s^{1/2} I'_r \Lambda_s^{1/2} \Pi_s^\star H_s$ and, if $\|H_s\| = 1$, we also have $\|\Pi_s^\star H_s\| = 1$ and thus this integrand is not bigger than $\lambda(r+1)_s$. Therefore, $\Delta(X, X')_t \leq L(r+1)_t$. Furthermore, $\sigma'_s \sigma'^\star_s = \Pi_s \Lambda_s I_r \Pi_s^\star$ is of rank $\leq r$, so $\sigma' \in \mathcal{S}_r$.

Now, let $\sigma'$ be any process in $\mathcal{S}_r$ and put $c'_s = \sigma'_s \sigma'^\star_s$. The kernel $K'_s$ of the linear map on $\mathbf{R}^d$ associated with the matrix $c'_s$ is of dimension at least $d - r$. The subspace $K_s$ of $\mathbf{R}^d$ generated by all eigenvectors of this linear map, which are associated with the eigenvalues $\lambda(1)_s, \ldots, \lambda(r+1)_s$, is of dimension at least $r+1$ (it is strictly bigger than $r+1$ if $\lambda(r+1)_s = \lambda(r+2)_s$). $K_s \cap K'_s$ is then not reduced to $\{0\}$ and we can thus find a process $H = (H_s)_{s \geq 0}$ with $\|H_s\| = 1$ and $H_s \in K_s \cap K'_s$ identically, and obviously this process can be chosen to be progressively measurable. Then $c'_s H_s = 0$ (because $H_s \in K'_s$) and $H_s^\star c_s H_s \geq \lambda(r+1)_s \|H_s\| = \lambda(r+1)_s$ (because $H_s \in K_s$). The first property above yields

$$H^\star \bullet X'_t = \int_0^t H_s^\star a_s \, \mathrm{d}s,$$

hence

$$\langle H^\star \bullet (X - X') \rangle_t = \langle H^\star \bullet X \rangle_t = \int_0^t (H_s^\star c_s^\star H_s) \, \mathrm{d}s,$$

which, by the second property above, is not less than $L(r+1)_t$. Hence, we are done. $\square$

In particular, $L(1)_t = \Delta(0; X)_t$ measures the 'distance' between the process $X$ and the pure drift process $X_0 + \int_0^t a_s \, \mathrm{d}s$. The following is obvious, with the convention $L(d+1)_t = 0$:

$$R_t \leq r \iff L(r+1)_t = \Delta(r; X)_t = 0, \tag{7}$$

where, similarly to (2), we have set, for all $t > 0$,

$$R(\omega)_t = \sup_{s \in [0,t]} \text{rank}(c_s(\omega)). \tag{8}$$

Hence, the random value $L(r+1)_t$ measures the distance between $X$ and the set of all processes with Brownian dimension $r$, over the time interval $[0, t]$, and *for our particular path* of $X$. Note also that it is an 'absolute' measure of the distance, which is multiplied by $u^2$ if we multiply the process $X$ by $u$.

Unfortunately, the variables $L(r)_t$ do not seem to be easy to 'estimate' from discrete observations since they involve eigenvalues. Hence, we will construct estimators below which are easier to handle.

## 3. Linear algebra preliminaries

Consider for a moment the toy model $X = \sigma W$, where $\sigma$ is a non-random $d \times q$ matrix. That is, we have (1) with $X_0 = 0$, $a_s = 0$ and $\sigma_s(\omega) = \sigma$, or, equivalently, $X$ is a Wiener process with covariance $\Sigma t$ at time t, where $\Sigma = \sigma \sigma^\star$. The observation scheme amounts to observing $n$ i.i.d. random vectors $G_i$, all of them $\mathcal{N}(0, \Sigma)$-distributed (viz. $G_i = \sqrt{n/T} \Delta_i^n X$, with the notation $\Delta_i^n X = X_{iT/n} - X_{(i-1)T/n}$). To infer $\text{rank}(\Sigma)$ from the observations of the $n$ first variables $G_k$, we can use the empirical covariance

$$\widehat{\Sigma}_n = \frac{1}{n} \sum_{k=1}^{n} G_k G_k^\star. \tag{9}$$

Indeed, the variables $G_k$ have a density over their support, which is a linear subspace with dimension $\text{rank}(\Sigma)$. Hence, $\text{rank}(\widehat{\Sigma}_n)$ is almost surely equal to $n$ when $n < \text{rank}(\Sigma)$, and to $\text{rank}(\Sigma)$ otherwise, and the problem is solved in a trivial way.

If we are in the same setting, except where $c_s$ depends on $s$ (and is still deterministic) and has a constant rank $r$, then, typically, the eigenspaces 'rotate' when $s$ varies, and the rank of $\widehat{\Sigma}_n$ above is a.s. equal to $d$ as soon as $n \geq d$. Therefore, $\widehat{\Sigma}_n$ gives no insight on the rank. So, the problem for non-homogeneous Wiener processes and, a fortiori, for general diffusions like (1), is actually more complex.

Despite the uselessness of the toy model consisting of an homogeneous Wiener process, let us give a couple of formulas relating to it, for further reference. We denote by $\mathcal{A}_r$ the family of all subsets of $\{1, \ldots, d\}$ with $r$ elements $(r = 1, \ldots, d)$. If $K \in \mathcal{A}_r$ and $\Sigma = (\Sigma^{ij}) \in \mathcal{M}_d$, we denote by $\det_K(\Sigma)$ the determinant of the $r \times r$ submatrix $(\Sigma^{kl} : k, l \in K)$ and let

$$\det(r; \Sigma) = \sum_{K \in \mathcal{A}_r} \det_K(\Sigma). \tag{10}$$

Observe that $\det(d; \Sigma) = \det(\Sigma)$, while $\det(1; \Sigma)$ is the trace of $\Sigma$.

**Lemma 3.** *If $\Sigma \in \mathcal{M}_d$ has eigenvalues $\lambda(1) \geq \cdots \geq \lambda(d) \geq 0$, we have, for $r = 1, \ldots, d$,*

$$\frac{1}{d(d-1)\cdots(d-r+1)} \det(r; \Sigma) \leq \lambda(1)\lambda(2)\cdots\lambda(r) \leq \det(r; \Sigma). \tag{11}$$

Notice that both inequalities in (11) may be equalities. It follows from (11) that, with the convention $0/0 = 0$, we have

$$1 \leq r \leq d \implies \begin{cases} r \leq \operatorname{rank}(\Sigma) & \implies & \det(r; \Sigma) > 0, \\ r > \operatorname{rank}(\Sigma) & \implies & \det(r; \Sigma) = 0, \end{cases} \tag{12}$$

$$2 \leq r \leq d \implies \frac{r!}{d!} \frac{\det(r; \Sigma)}{\det(r-1; \Sigma)} \leq \lambda(r) \leq \frac{d!}{(r-1)!} \frac{\det(r; \Sigma)}{\det(r-1; \Sigma)}. \tag{13}$$

**Proof of Lemma 3.** We expand the characteristic polynomial of $\Sigma$ as

$$\det(\Sigma - \lambda I) = (-\lambda)^d + \sum_{r=1}^{d} (-\lambda)^{d-r} \det(r; \Sigma).$$

In view of the well-known expressions for the 'symmetrical functions' of the roots of a polynomial, we get

$$\sum_{1 \leq i_1 < \cdots < i_r \leq d} \lambda(i_1)\lambda(i_2)\cdots\lambda(i_r) = \det(r; \Sigma)$$

and thus both sides of (11) are obvious. □

Next, we consider a sequence $(G_i)_{i \geq 1}$ of i.i.d. $\mathcal{N}(0, \Sigma)$-distributed random vectors. For all $j = 1, \ldots, d$, we define two random elements of $\mathcal{M}_d$ by

$$\zeta_j = \sum_{i=1}^{j} G_i G_i^\star, \qquad \zeta_j' = \sum_{i=d+1}^{d+j} G_i G_i^\star \tag{14}$$

and we consider the mean and covariance of the random vector $(\det(r; \zeta_r)/r! : 1 \leq r \leq d)$:

$$\left. \begin{array}{l} \gamma(r; \Sigma) = \dfrac{1}{r!} \mathbf{E}(\det(r; \zeta_r)), \\ \Gamma(r, r'; \Sigma) = \dfrac{1}{r!r'!} \mathbf{E}(\det(r; \zeta_r) \det(r'; \zeta_{r'})) - \gamma(r; \Sigma)\gamma(r'; \Sigma). \end{array} \right\} \tag{15}$$

Since $(\zeta_j, \zeta_j' : 1 \leq j \leq d)$ are i.i.d., we also have

$$\Gamma(r, r'; \Sigma) = \frac{1}{r!r'!} \mathbf{E}(\det(r; \zeta_r) \det(r'; \zeta_{r'}) - \det(r; \zeta_r) \det(r'; \zeta_{r'}')). \tag{16}$$

**Lemma 4.** *If $r \in \{1, \ldots, d\}$, we have*

$$\gamma(r; \Sigma) = \det(r; \Sigma). \tag{17}$$

*Moreover,*

$$\left.\begin{array}{rl} r \leq \operatorname{rank}(\Sigma) & \implies \gamma(r; \Sigma) > 0, \Gamma(r, r; \Sigma) > 0, \\ r > \operatorname{rank}(\Sigma) & \implies \gamma(r; \Sigma) = 0, \Gamma(r, r; \Sigma) = 0. \end{array}\right\} \tag{18}$$

**Proof.** To prove (17), it is enough to show that for any $K \in \mathcal{A}_r$, we have $\mathbf{E}(\det_K(\zeta_r)) = r! \det_K(\Sigma)$, and for this, it is no restriction to assume that $K = \{1, \ldots, r\}$. We denote by $\mathcal{P}_r$ the set of all permutations of the set $\{1, \ldots, r\}$ and by $\varepsilon(\tau)$ the signature of the permutation $\tau$. Then

$$\det_K(\zeta_r) = \sum_{\tau \in \mathcal{P}_r} (-1)^{\varepsilon(\tau)} \prod_{l=1}^r \zeta_r^{l\tau(l)} = \sum_{1 \leq k_1, \ldots, k_r \leq r} \sum_{\tau \in \mathcal{P}_r} (-1)^{\varepsilon(\tau)} \prod_{l=1}^r G_{k_l}^l G_{k_l}^{\tau(l)}$$

and each summand of the first sum of the extreme right-hand side is the determinant of a matrix with rank less than $r$, unless all $k_l$ are distinct. So, it is enough to sum over all $r$-tuples $(k_1, \ldots, k_r)$ with distinct entries between 1 and $r$, that is, for $r$-tuples with $k_i = \tau'(i)$ for some $\tau' \in \mathcal{P}_r$. In other words, we have

$$\det_K(\zeta_r) = \sum_{\tau' \in \mathcal{P}_r} \sum_{\tau \in \mathcal{P}_r} (-1)^{\varepsilon(\tau)} \prod_{l=1}^r G_{\tau'(l)}^l G_{\tau'(l)}^{\tau(l)}. \tag{19}$$

Since the variables $G_n$ are independent and $\mathbf{E}(G_n^k G_n^l) = \Sigma^{kl}$, we deduce that

$$\mathbf{E}(\det_K(\zeta_r)) = r! \sum_{\tau \in \mathcal{P}_r} (-1)^{\varepsilon(\tau)} \prod_{l=1}^r \Sigma^{l\tau(l)} = r! \det_K(\Sigma)$$

and we have (17).

If $r > \operatorname{rank}(\Sigma)$, we have $\det(r; \Sigma) = 0$ (see (12)): the non-negative variable $\det(r; \zeta_r)$ has zero expectation, so it is a.s., null and we have the second part of (18). Finally, let $r \leq \operatorname{rank}(\Sigma)$. Again, by (12), we have $\mathbf{E}(\det(r; \zeta_r)) > 0$. Also, observe that $\det(r; \zeta_r)$ is a continuous function of the random vectors $G_n$ for $n = 1, \ldots, r$ which vanishes if all these $G_n$ are 0. Thus, $\det(r; \zeta_r)$ can take arbitrarily small values with positive probability and it has a positive expectation, so it is not degenerate and we get the first part of (18). □

## 4. Limit theorems for estimators of the Brownian dimension

It turns out that determinants or 'integrated determinants' are much easier to estimate than eigenvalues or integrated eigenvalues. So, in view of (11) and (12), one might replace the variable $L(r)_t$ of (6) by

$$L'(r)_t = \begin{cases} \int_0^t \det(1; c_s) \, ds, & \text{if } r = 1, \\ \int_0^t \frac{\det(r; c_s)}{\det(r-1; c_s)} \, ds, & \text{if } r \geq 2. \end{cases}$$

However, $L'(r)_t$ for $r \geq 2$ is still not so easy to estimate. For example, for the toy model of Section 3, the variable $\det(r; \zeta_r)/r!$ is an unbiased estimator of $\det(r; \Sigma)$ (see (17)), but we have no explicit unbiased estimator for a quotient like $\det(r; \Sigma)/\det(r-1; \Sigma)$.

So, we propose to measure *the distance between $X$ and the set of models with multiplicity $r$*, over the time interval $[0, t]$, by the following random variable:

$$\overline{L}(r)_t = \int_0^t \det(r; c_s) \, ds. \tag{20}$$

Up to multiplicative constants, this more or less amounts to replacing the 'natural' distance $L(r)_t$ by $\int_0^t \lambda(1)_s \cdots \lambda(r)_s \, ds$. The variables $L(r)_t$, $L'(r)_t$ and $\overline{L}(r)_t$ convey essentially the same information as far as rank is concerned and, in particular, they vanish simultaneously, which is the most important property for our purposes. In other words, exactly as in (7), we have

$$R_t \leq r \iff \overline{L}(r+1)_t = 0. \tag{21}$$

By virtue of (17), we can rewrite $\overline{L}(r)$ as follows (we also introduce additional variables $Z(r, r')$, using the notation (15)):

$$\overline{L}(r)_t = \int_0^t \gamma(r; c_s) \, ds, \qquad Z(r, r')_t = \int_0^t \Gamma(r, r'; c_s) \, ds. \tag{22}$$

Now, we need to approximate the variables in (22) by variables which depend on our discrete observations only. To this end, we introduce the random matrices

$$\zeta(r)_i^n = \sum_{j=1}^r (\Delta_{i+j-1}^n X)(\Delta_{i+j-1}^n X)^\star, \qquad \text{where } \Delta_i^n X = X_{iT/n} - X_{(i-1)T/n}. \tag{23}$$

We have $\zeta(r)_i^n \in \mathcal{M}_d$ and $\text{rank}(\zeta(r)_i^n) \leq r$. We then set (with $[x]$ being the integer part of $x$)

$$\overline{L}(r)_t^n = \frac{n^{r-1}}{T^{r-1} r!} \sum_{i=1}^{[nt/T]-r+1} \det(r; \zeta(r)_i^n) \tag{24}$$

and

$$Z(r, r')_t^n = \frac{n^{r+r'-1}}{T^{r+r'-1} r! r'!} \sum_{i=1}^{[nt/T]-d-r'+1} (\det(r; \zeta(r)_i^n) \det(r'; \zeta(r')_i^n)$$

$$- \det(r; \zeta(r)_i^n) \det(r'; \zeta(r')_{d+i}^n)). \tag{25}$$

The first key theorem concerns the 'consistency' of these variables.

**Theorem 5.** *Under* (H), *the variables* $\overline{L}(r)_t^n$ *and* $Z(r, r')_t^n$ *converge in probability to* $\overline{L}(r)_t$ *and* $Z(r, r')_t$, *respectively, uniformly in* $t \in [0, T]$.

This is not enough for our purposes and we need rates of convergence. For this, (H) is not sufficient and some additional regularity on the coefficients $a$ and $\sigma$ is necessary. A first set of sufficient conditions is simple enough.

*Hypothesis (H1).* We have (H) with a càdlàg process $a$ and a process $\sigma$ which is Hölder continuous (in time) with index $\rho > 1/2$, in the sense that

$$\sup_{0 \leq s < t \leq T} \frac{\|\sigma_t - \sigma_s\|}{(t-s)^\rho} < \infty \qquad \text{a.s.} \tag{26}$$

The above assumption on $a$ is quite mild and the assumption on $\sigma$ is reasonable when $\sigma$ is deterministic. However, in the pure diffusion case, we have $\sigma_s = \sigma(X_s)$ for, say, a Lipschitz or locally Lipschitz function $\sigma$ and, of course, (26) fails for any $\rho \geq 1/2$. This assumption also fails when $\sigma_s$ is a 'stochastic volatility' driven by an Itô equation and even more if this equation has jumps!

Therefore, for practical purposes which are especially relevant in finance, we need to replace (H1) by a different assumption. This assumption looks (is ?) complicated to state, but it essentially says that $a$ is as in (H1) and that the process $\sigma$ follows a jump-diffusion Itô equation, or, in other words, that it is driven by a Wiener process and a Poisson random measure). In particular, it is satisfied in the pure diffusion case when $\sigma_s = \sigma(X_s)$, with a $C^2$ function $\sigma$.

*Hypothesis (H2).* We have (H), the process $a$ is càdlàg and the process $\sigma$ is a (possibly discontinuous) Itô semimartingale on $[0, T]$, that is, for $t \leq T$, we have

$$\begin{aligned}
\sigma_t = \sigma_0 &+ \int_0^t a'_s \, ds + \int_0^t \sigma'_{s-} \, dW_s \\
&+ \int_0^t \int_E \varphi \circ w(s-, x)(\mu - \nu)(ds, dx) + \int_0^t \int_E (w - \varphi \circ w)(s-, x) \mu(ds, dx).
\end{aligned} \tag{27}$$

Here, $\sigma'$ is $\mathbf{R}^d \otimes \mathbf{R}^q \otimes \mathbf{R}^q$-valued adapted càdlàg and $a'$ is $\mathbf{R}^d \otimes \mathbf{R}^q$-valued predictable and locally bounded; $\mu$ is a Poisson random measure on $(0, \infty) \times E$, independent of $W$ and $V$, with intensity measure $\nu(dt, dx) = dt F(dx)$ with $F$ a $\sigma$-finite measure on some Polish space $(E, \mathcal{E})$; $\varphi$ is a continuous function on $\mathbf{R}^{dq}$ with compact support, which coincides with the identity on a neighborhood of 0; finally, $w(\omega, s, x)$ is a map $\Omega \times [0, \infty) \times E \to \mathbf{R}^d \otimes \mathbf{R}^q$ which is $\mathcal{F}_s \otimes \mathcal{E}$-measurable in $(\omega, x)$ for all $s$, cadlag in $s$ and such that $\int_E (1 \wedge \sup_{\omega \in \Omega, s \leq T} \|w(\omega, s, x)\|^2) F(dx) < \infty$.

These conditions are indeed quite easy to check in practice. They accommodate the case of a stochastic volatility driven by a Wiener process having some (or all) components independent of $X$: since $W$ has an 'arbitrary' dimension $q$ in this paper, possibly $q > d$, there might be components used for $X$ in (1) and other components used in (27).

**Theorem 6.** *Assume either* (H1) *or* (H2). *The $d$-dimensional processes $(V(r)^n_t)_{1 \leq r \leq d}$ with components*

$$V(r)^n_t = \sqrt{n}(\overline{L}(r)^n_t - \overline{L}(r)_t) \tag{28}$$

*converge stably in law to a limiting process $(V(r)_t)_{1 \leq r \leq d}$ which is defined on an extension of the original space and which, conditionally on $\mathcal{F}$, is a non-homogeneous Wiener process with quadratic variation process $t \mapsto (TZ(r,r')_t)_{1 \leq r, r' \leq d}$.*

**Proof of Theorems 5 and 6.** The proof involves several steps.

(1) It is based on the following two results of [4]. Take $N$ functions $g_j$ on $\mathbf{R}^d$ which are $C^2$ and with polynomial growth and even. Set

$$Y(g_1, \ldots, g_N)^n_t = \frac{T}{n} \sum_{i=1}^{[nt/T]-N+1} \prod_{k=1}^{N} g_k(\sqrt{n/T}\Delta^n_{i+k-1}X).$$

Then, under (H), we have

$$Y(g_1, \ldots, g_N)^n_t \xrightarrow{\mathbf{P}} Y(g_1, \ldots, g_N)_t := \int_0^t y(g_1, \ldots, g_N; c_s)\,\mathrm{d}s,$$

where the convergence in uniform in $t \in [0, T]$ and where $y(g_1, \ldots, g_N; \Sigma)$ is, for any $d \times d$ covariance matrix $\Sigma$, the expectation of the variable

$$\gamma(g_1, \ldots, g_N) = \prod_{k=1}^{N} g_k(G_k)$$

and the $G_n$'s are i.i.d. random vectors with law $\mathcal{N}(0, \Sigma)$, as in Section 3.

If, further, (H2) holds, then for any array $((g^j_1, \ldots, g^j_{N_j}) : 1 \leq j \leq J)$ with $g^j_i$ as above, the $J$-dimensional processes $(\sqrt{n/T}(Y(g^j_1, \ldots, g^j_{N_j})^n_t - Y(g^j_1, \ldots, g^j_{N_j})_t))_{1 \leq j \leq J}$ converge stably in law to a limiting process which, conditionally on $\mathcal{F}$, is a non-homogeneous Wiener process with quadratic variation process $\int_0^t \Gamma(c_s)\,\mathrm{d}s$ and where $\Gamma(\Sigma)$ is the covariance matrix of the random vector $(\gamma(g^j_1, \ldots, g^j_{N_j}))_{1 \leq j \leq J}$, as defined above. Under (H1) instead of (H2), the same result holds: it is not explicitly stated in [4], but the proof is similar, and technically much simpler.

(2) These results extend by 'linearity' in an obvious way. More precisely, for $1 \leq j \leq J$, set

$$Y(j)^n_t = \frac{T}{n} \sum_{i=1}^{[nt/T]-N_j+1} h_j(\sqrt{n/T}\Delta^n_i X, \sqrt{n/T}\Delta^n_{i+1}X, \ldots, \sqrt{n/T}\Delta^n_{i+N_j-1}X), \tag{29}$$

where each $h_j$ is a linear combination of tensor products $g_1 \otimes \cdots \otimes g_{N_j}$, where the $g_i$'s are $C^2$ functions on $\mathbf{R}^d$, even and with polynomial growth. Also, denote by $M(\Sigma)$ and $C(\Sigma)$ the mean vector and the covariance matrix, respectively, of the $J$-dimensional random vector $(h_j(G_1, \ldots, G_{N_j}))_{1 \leq j \leq J}$, with $G_i$ as above. Then:

1. under (H), we have

$$Y(j)_t^n \xrightarrow{\mathbf{P}} Y(j)_t := \int_0^t M^j(c_s)\,\mathrm{d}s, \qquad \text{uniformly in } t \in [0,T]; \tag{30}$$

2. under (H1) or (H2), the $J$-dimensional processes with components $\sqrt{n/T}(Y(j)_t^n - Y(j)_t)$ converge stably in law to a limiting process which, conditionally on $\mathcal{F}$, is a non-homogeneous Wiener process with quadratic variation process $\int_0^t C(c_s)\,\mathrm{d}s$.

(3) The theorem is now almost trivial. The determinants entering (24) (resp., (25)) are sums of even monomials of the components of $\Delta_{i+j-1}^n X$ for $1 \leq j \leq 2d$, each one with degree $2r$ (resp., $2(r+r')$). More specifically, $\overline{L}(r)^n$ is of type (29), with $N_r = r$ and the function

$$h_r(x_1, \ldots, x_r) = \frac{1}{r!} \det\left(r, \sum_{j=1}^r x_j x_j^\star\right),$$

whereas $Z(r,r')^n$ is of type (29), with $N_{r,r'} = d + r'$ and the function

$$h_{r,r'}(x_1, \ldots, x_{d+r'}) = \frac{1}{r!r'!}\left(\det\left(r, \sum_{j=1}^r x_j x_j^\star\right) \det\left(r', \sum_{j=1}^{r'} x_j x_j^\star\right) - \det\left(r, \sum_{j=1}^r x_j x_j^\star\right) \det\left(r', \sum_{j=d+1}^{d+r'} x_j x_j^\star\right)\right).$$

Theorem 5 then follows readily from step 2 and the relations (15). □

In the sequel, we will also need some estimates on the moments of $V(r)_t^n$, uniform in $n$. It follows from the proofs in [4] that, under (H2) and for each $t \in (0,T]$, there is a sequence $A_{p,t}$ of $\mathcal{F}_t$-measurable sets such that

$$\left.\begin{array}{l} A_{p,t} \uparrow \Omega \text{ as } p \to \infty, \\ p \geq 1,\ n \geq 1,\ t \in (0,T] \quad \Longrightarrow \quad \mathbf{E}(|V(r)_t^n|^2 \mathbf{1}_{A_{p,t}}) \leq C_p \end{array}\right\} \tag{31}$$

for a suitable sequence of constants $C_p$ (depending on $T$). The same result also holds under (H1).

Now, if the coefficient $a$ is bounded and if we have (H1) with (26) holding uniformly in $\omega$, or (H2) with $a_s'(\omega)$, $\sigma_s' - \omega$ bounded and $\|w(\omega,s,x)\| \leq h(x)$ for some function having $\int_E (1 \wedge h(x)^2) F(\mathrm{d}x) < \infty$, one can (easily) prove that (31) holds for $A_{1,t} = \Omega$, and so there is a constant $C$, depending on $T$, such that

$$n \geq 1, \ t \in (0, T] \quad \Longrightarrow \quad \mathbf{E}(|V(r)_t^n|^2) \leq C. \tag{32}$$

## 5. Tests based on thresholds

### 5.1. A test based on an absolute threshold

We return to the initial problem, in the light of the second comment of Section 1. Namely, we want to decide which integer value (between 0 and $d$) the variable $R$ of (2) takes for the particular path $\omega$ which is known only through the observations $X_{iT/n}$. In principle, we have our observations $X_{iT/n}$ for $i = 0, \ldots, n$, but it may be interesting to determine how our estimators behave as time changes. This is why we also give estimators for the variable $R_t$ of (8), based on the observation of $X_{iT/n}$ for $i = 0, \ldots, [nt/T]$.

Let us reiterate that in the ideal scheme (the whole path of $X$ is known over $[0, T]$), we also know $R = R(\omega)$, whereas we have the equivalence (21). In view of this, and taking into account the convergence result in Theorem 5, it seems natural to operate as follows. We choose a sequence of positive numbers $\rho_n$ such that

$$\rho_n \to 0, \qquad \rho_n \sqrt{n} \to \infty. \tag{33}$$

We then take the following 'estimator' for $R_t$:

$$\widehat{R}_{n,t} = \inf(r \in \{0, \ldots, d-1\} : \overline{L}(r+1)_t^n < \rho_n t), \tag{34}$$

with $\inf(\varnothing) = d$. That this estimator is a priori reasonable is due to the fact that if we set $R_{n,t} = \inf(r \in \{0, \ldots, d-1\} : \overline{L}(r+1)_t < \rho_n t)$, then by (21) and the property $\rho_n \to 0$, we have $\mathbf{P}(R_{n,t} = R_t) \to 1$ as $n \to \infty$. We take a threshold of the form $\rho_n t$ because $\overline{L}(r+1)_t^n$ is roughly proportional to $t$.

**Remark 7.** Another equally reasonable estimator, which is a kind of 'dual' of $\widehat{R}_{n,t}$, is the following one:

$$\widehat{R}'_{n,t} = \sup(r \in \{1, \ldots, d\} : \overline{L}(r)_t^n \geq \rho_n t), \tag{35}$$

with $\sup(\varnothing) = 0$. The analysis of $\widehat{R}_{n,t}$ below carries over for $\widehat{R}'_{n,t}$ in a very similar way.

**Remark 8.** The choice of the threshold $\rho_n$ is arbitrary, upon the fact that (33) holds; asymptotically all choices are equivalent. In practice, though, it is of primary importance because $n$, albeit large, is given and is of course not infinite! Even worse, an absolute threshold like in (34) is sensitive to the unit in which the values of the $X_t^i$ are expressed. For example, if we multiply all components by the same (known) constant, the estimator of the Brownian dimension provides a different value. So, using an absolute threshold is probably *not* advisable in general. Nevertheless, we pursue here the analysis of tests based on an absolute threshold since they may serve as a case study and are somewhat simpler to study than the tests based on relative thresholds which are introduced later.

The integer-valued estimator $\widehat{R}_{n,t}$ should be analyzed using the testing methodology rather than as a usual estimator: we test the hypothesis $R_t = r$ with the critical region $\{\widehat{R}_{n,t} \neq r\}$. The 'power function' is, in principle, the probability of rejection, a function of the underlying probability measure. Here, we have a single $\mathbf{P}$ and $R_t$ is (possibly) random. We thus develop two different substitutes for the power function.

## 5.2. A first substitute for the power function

A seemingly acceptable version of the power function is

$$\widehat{\beta}_{n,t}^r(r') = \mathbf{P}(\widehat{R}_{n,t} \neq r \mid R_t = r'), \qquad r' = 0, 1, \ldots, d, \tag{36}$$

provided $\mathbf{P}(R_t = r') > 0$. We explicitly mention the number $n$ of observations and the number $r$, but it also depends on the sequence $\rho_n$. The index $r$ indicates the 'test' with null hypothesis $R_t = r$ which we are performing, while the index $r'$ indicates the 'true' value or $R_t$. So, $\widehat{\beta}_{n,t}^r(r)$ should be small and $\widehat{\beta}_{n,t}^r(r')$ should be close to 1 when $r' \neq r$.

We have a first – quite simple – result.

**Theorem 9.** *Under (33) and either* (H1) *or* (H2), *we have, for all $r, r'$ in $\{1, \ldots, d\}$, and provided $\mathbf{P}(R_t = r') > 0$,*

$$\widehat{\beta}_{n,t}^r(r') \longrightarrow \begin{cases} 1, & \text{if } r \neq r', \\ 0, & \text{if } r = r'. \end{cases} \tag{37}$$

Another equivalent (simpler) way of stating this result is to write

$$\mathbf{P}(\widehat{R}_{n,t} \neq R_t) \to 0. \tag{38}$$

This is more intuitive, but somehow further away from the way results on tests are usually stated.

**Proof of Theorem 9.** For each $s = 1, \ldots, (r'+1) \wedge d$, we set $\delta_{n,t}^s(r') = \mathbf{P}(\overline{L}(s)_t^n < \rho_n, R_t = r')$. Observe that, if $\rho_n' = \rho_n \sqrt{n}$, and with the notation (28),

$$\left.\begin{array}{l}\delta_{n,t}^s(r') \leq \mathbf{P}(\overline{L}(s)_t < 2\rho_n t, R_t = r') + \mathbf{P}(|V(s)_t^n| > \rho_n' t), \\ \mathbf{P}(R_t = r') - \delta_{n,t}^s(r') \leq \mathbf{P}(\overline{L}(s)_t > \rho_n t/2, R_t = r') + \mathbf{P}(|V(s)_t^n| \geq \rho_n' t/2).\end{array}\right\} \tag{39}$$

Theorem 6 yields that the sequence $V(s)_t^n$ converges in law when $n$ goes to infinity, whereas on the set $\{R_t = r'\}$, we have $\overline{L}(s)_t > 0$ if $s \leq r'$, and $\overline{L}(s)_t = 0$ if $s > r'$. Therefore, it follows from (33) that

$$\left.\begin{array}{l}s \leq r' \implies \delta_{n,t}^s(r') \to 0, \\ s = r' + 1 \leq d \implies \delta_{n,t}^s(r') \to \mathbf{P}(R_t = r').\end{array}\right\} \tag{40}$$

Now, $\{\widehat{R}_{n,t} \neq r'\} = (\bigcup_{1 \leq s \leq r'} \{\overline{L}(s)^n_t < \rho_n t\}) \cup \{\overline{L}(r'+1)^n_t \geq \rho_n t\}$, with the convention that $\{\overline{L}(d+1)^n_t \geq \rho_n t\} = \varnothing$. Then

$$\mathbf{P}(\widehat{R}_{n,t} \neq r' = R_t) \leq \begin{cases} \displaystyle\sum_{s=1}^{r'} \delta^s_{n,t}(r') + \mathbf{P}(R_t = r') - \delta^{r'+1}_{n,t}(r'), & \text{if } r' \leq d-1, \\ \displaystyle\sum_{s=1}^{r'} \delta^s_{n,t}(r'), & \text{if } r' = d. \end{cases}$$

It then readily follows from (40) that $\widehat{\beta}^{r'}_{n,t}(r') \to 0$ as soon as $\mathbf{P}(R_t = r') > 0$. Under this assumption, and if $r \neq r'$, we clearly have $\widehat{\beta}^r_{n,t}(r') = 1 - \mathbf{P}(\widehat{R}_{n,t} = r \mid R_t = r') \geq 1 - \widehat{\beta}^{r'}_{n,t}(r)$, hence $\widehat{\beta}^r_{n,t}(r') \to 1$. $\square$

The previous result seems to settle the matter. However, it is not as nice as it may look, because it gives no 'rate' for the convergence in (37) or (38) and is thus impossible to use in practice. The impossibility of getting a rate is apparent in (39): the second terms on the right may be more or less controlled through estimates like (31), but the first terms on the right cannot be controlled at all; indeed, if $R_t = r$, the variable $\overline{L}(r)_t$ is positive, but may be arbitrarily close to 0.

### 5.3. A second substitute for the power function

As emphasized in Comment 4 of Section 1, one may reasonably decide that $R_t = r$ if $r$ is the 'true' Brownian dimension in a 'significant' way, which means, in particular, that the 'distance' between the model $X$ and the set of models with Brownian dimension $r' < r$ is not 'infinitesimal'. This may be interpreted as the property that $\overline{L}(r)_t$ exceeds some positive level for all $r \leq R_t$.

In other words, we set $B_{r',\varepsilon,t} = \{R_t = r', \overline{L}(r)_t \geq \varepsilon t \text{ for } r = 1, \ldots, r'\}$ and define the 'power function' as being

$$\widehat{\beta}^r_{n,t}(r', \varepsilon) = \mathbf{P}(\widehat{R}_{n,t} \neq r \mid B_{r',\varepsilon,t}), \tag{41}$$

provided $\mathbf{P}(B_{r',\varepsilon,t}) > 0$.

Evaluating $\beta^r_{n,t}(r', \varepsilon)$ is still difficult, because it involves the unknown quantity $\mathbf{P}(B_{r',\varepsilon,t})$. So, we provide a result which does not directly give the power function itself, but which is probably more relevant for applications.

**Theorem 10.** *Under* (H1) *or* (H2), *there are $\mathcal{F}_t$-measurable sets $(A_{p,t})$ increasing to $\Omega$ as $p \to \infty$ and constants $C_p$ such that, for all $r$ in $\{1, \ldots, d\}$, and provided $\rho_n < \varepsilon/2$,*

$$\mathbf{P}(\{\widehat{R}_{n,t} \neq r\} \cap B_{r,\varepsilon,t} \cap A_{p,t}) \leq \frac{C_p}{n\rho_n^2} \tag{42}$$

*for all $t \in (0, T]$. If, further, (32) holds, we can find a constant such that*

$$\mathbf{P}(\{\widehat{R}_{n,t} \neq r\} \cap B_{r,\varepsilon,t}) \leq \frac{C}{n\rho_n^2} \tag{43}$$

*or, in other words, the 'level' satisfies $\widehat{\beta}_{n,t}^r(r, \varepsilon) \leq C/(n\rho_n^2 \mathbf{P}(B_{r,\varepsilon,t}))$.*

Note that we can choose $\rho_n$ above at will, provided it satisfies (33), and none of $A_p$, $B_{r',\varepsilon,t}$, $C_p$ depend on this choice (the estimator $\widehat{R}_{n,t}$ does, though). So, we can obtain a rate $1/n^\theta$ for any $\theta \in (0, 1)$, as close to 1 as one wishes.

**Proof of Theorem 10.** We consider the sets $A_{p,t}$ for which (31) holds and denote by $C'_p$ the constants occurring in that formula. For $s = 1, \ldots, (r+1) \wedge d$, we set ($\varepsilon > 0$ being fixed) $\delta_{n,t}^s(r, p) = \mathbf{P}(\{\overline{L}(s)_t^n < \rho_n t\} \cap B_{r,\varepsilon,t} \cap A_{p,t})$. Exactly as for (39), we have

$$\delta_{n,t}^s(r, p) \leq \mathbf{P}(\{\overline{L}(s)_t < 2\rho_n t\} \cap B_{r,\varepsilon,t}) + \mathbf{P}(\{|V(s)_t^n| > \rho_n' t\} \cap A_{p,t}),$$

$$\mathbf{P}(B_{r,\varepsilon,t} \cap A_{p,t}) - \delta_{n,t}^{\prime s}(r, p) \leq \mathbf{P}(\{\overline{L}(s)_t > \rho_n t/2\} \cap B_{r,\varepsilon,t}) + \mathbf{P}(\{|V(s)_t^n| \geq \rho_n' t/2\} \cap A_{p,t}).$$

Taking into account (31), $\rho_n < \varepsilon/2$ and the facts that $\overline{L}(s)_t = 0$ if $R_t = r < s$ and that $\overline{L}(r)_t \geq \varepsilon t$ on $B(r, \varepsilon, t)$, we deduce from Chebyshev's inequality that

$$s \leq r \quad \Longrightarrow \quad \delta_{n,t}^s(r, p) \leq \frac{C'_p}{\rho_n'^2}, \qquad \mathbf{P}(B_{r,\varepsilon,t} \cap A_{p,t}) - \delta_{n,t}^{r+1}(r, p) \leq \frac{4C'_p}{\rho_n'^2}, \tag{44}$$

where the second equality makes sense when $r < d$ only. Once more applying the identity $\{\widehat{R}_{n,t} \neq r\} = \{\overline{L}(r+1)_t^n \geq \rho_n t\} \cup (\bigcup_{1 \leq s \leq r} \{\overline{L}(s)_t^n < \rho_n t\})$, we get

$$\mathbf{P}(\{\widehat{R}_{n,t} \neq r\} \cap B_{r,\varepsilon,t} \cap A_{p,t}) \leq \begin{cases} \sum_{s=1}^r \delta_{n,t}^s(r, p) + \mathbf{P}(B_{r,\varepsilon,t} \cap A_{p,t}) - \delta_{n,t}^{r+1}(r, p), & \text{if } r < d, \\ \sum_{s=1}^r \delta_{n,t}^s(r, p), & \text{if } r = d. \end{cases}$$

We then deduce (42) from (44) if we put $C_p = (4 + d)C'_p$. Finally, under (32), we may choose $A_{1,t} = \Omega$ above and thus (43) with $C = (4 + d)C'_1$. $\square$

Of course, (42) is not useful in general, although it gives us a rate, because we do not know the sets $A_{p,t}$. If (32) holds, the result appears much more satisfactory; however, we still do not know the constant $C$ in (43) and have no mean to guess what it is from the observations.

### 5.4. Tests based on a relative threshold

In practice, the previous tests are not recommended; see Remark 8. We now exhibit other tests which are scale-invariant.

If we multiply $X$ by a constant $\delta > 0$, then $c_s$ is multiplied by $\delta^2$ and both $\overline{L}(r)_t^n$ and $\overline{L}(r)_t$ are multiplied by $\delta^{2r}$. Then, for any given sequence $\rho_n \in (0,1]$ satisfying (33), the following two 'estimators' of $R$, which are candidates to be *explicative Brownian dimensions*, are scale-invariant:

$$\left.\begin{array}{l}\widetilde{R}_{n,t} = \inf(r \in \{0\ldots, d-1\}: \overline{L}(r+1)_t^n < \rho_n t^{-1/r}(\overline{L}(r)_t^n)^{(r+1)/r}),\\ \widetilde{R}'_{n,t} = \inf(r \in \{0,\ldots, d-1\}: \overline{L}(r+1)_t^n < \rho_n t^{-r}(\overline{L}(1)_t^n)^{r+1}),\end{array}\right\} \quad (45)$$

with the convention that $\overline{L}(0)_t^n = 1$ and, again, $\inf(\varnothing) = d$. The presence of $t^{1/r}$ or $t^r$ above accounts for the fact that $\overline{L}(r)_t^n$ is roughly proportional to $t$, as in (34).

Note that $\widetilde{R}'_{n,t} \geq 1$, even when $R_t = 0$. So, if $R_t = 0$, this estimator is bad, but, in this case, our problem is essentially meaningless anyway! When $R_t \geq 1$, the significance of these two estimators is essentially as follows: $\widetilde{R}_{n,t}$ is the smallest integer $r$ for which there is a 'large' drop between the explicative powers of the models with Brownian dimensions $r$ and $r+1$, whereas $\widetilde{R}'_{n,t}$ is the smallest integer $r$ at which the ratio between the contributions of the $(r+1)$th and the first Brownian dimension is smaller than $\rho_n$. Clearly, there exist other estimators of the same kind, with slightly different meanings, the above two being the extremes. All such estimators are amenable to essentially the same mathematical analysis.

In practice, the choice of $\rho_n$ is relative to the physical phenomenon under consideration and to the use which is made of the model (prediction, simulation, computation of extreme values, etc.). Roughly speaking, the choice should reflect the physical effects which are modelled as the driving noise and the intensity of components of the noise which are considered important to capture essential properties of the model.

Here, again, the substitutes for the power functions are

$$\widetilde{\beta}_{n,t}^r(r') = \mathbf{P}(\widetilde{R}_{n,t} \neq r \mid R_t = r'), \qquad \widetilde{\beta}_{n,t}^{\prime r}(r') = \mathbf{P}(\widetilde{R}'_{n,t} \neq r \mid R_t = r'). \quad (46)$$

We now aim to produce a result similar to Theorem 9.

**Theorem 11.** *Under (33) and either* (H1) *or* (H2), *we have, for all* $r, r'$ *in* $\{1, \ldots, d\}$, *and provided* $\mathbf{P}(R_t = r') > 0$.

$$\widetilde{\beta}_{n,t}^r(r') \longrightarrow \begin{cases} 1, & \text{if } r \neq r', \\ 0, & \text{if } r = r' \end{cases} \quad (47)$$

*and the same for* $\widetilde{\beta}_{n,t}^{\prime r}(r')$.

**Proof.** We prove the result for $\widetilde{\beta}_{n,t}^r(r')$ only, the other case being similar. For each $s = 1, \ldots, r' \wedge d$, we set

$$\delta_{n,t}^s(r') = \mathbf{P}(\overline{L}(s)_t^n < \rho_n t^{-1/s}(\overline{L}(s-1)_t^n)^{s/(s-1)}, R_t = r').$$

As in Theorem 9, we write $\rho'_n = \rho_n\sqrt{n}$, and with the convention $0/0 = 1$ and $t \in (0,T]$ fixed, we put

$$\widetilde{V}(s)_n = t^{1/s}\sqrt{n}\left(\frac{\overline{L}(s)_t^n}{(\overline{L}(s-1)_t^n)^{s/(s-1)}} - \frac{\overline{L}(s)_t}{(\overline{L}(s-1)_t)^{s/(s-1)}}\right). \qquad (48)$$

Observe that, similarly to (39), we have

$$\delta_{n,t}^s(r') \leq \mathbf{P}(\overline{L}(s)_t < 2\rho_n t^{-1/s}(\overline{L}(s-1)_t)^{(s-1)/s}, R_t = r') + \mathbf{P}(|\widetilde{V}(s)_n| > \rho'_n, R_t = r'),$$

$$\mathbf{P}(R_t = r') - \delta_{n,t}^s(r') \leq \mathbf{P}(\overline{L}(s)_t > \rho_n t^{-1/s}(\overline{L}(s-1)_t)^{(s-1)/s}/2, R_t = r')$$
$$+ \mathbf{P}(|\widetilde{V}(s)_n| \geq \rho'_n/2, R_t = r').$$

Moreover, Theorem 6 yields that on the set $\{\overline{L}(s-1)_t > 0\} = \{R_t \geq s-1\}$, the variables $\widetilde{V}(s)_n$ converge stably in law to the variable

$$\widetilde{V}(s) = \frac{t^{1/s}\overline{L}(s)_t}{(\overline{L}(s-1)_t)^{s/(s-1)}}\left(\frac{V(s)_t}{\overline{L}(s)_t} - \frac{s}{s-1}\frac{V(s-1)_t}{\overline{L}(s-1)_t}\right).$$

Therefore, since (33) holds, we get (40). At this stage, we can reproduce the end of the proof of Theorem 9 to obtain (47). $\square$

Obviously, the comments made after Theorem 9 for the estimators $\widehat{R}_{n,t}$ apply to $\widetilde{R}_{n,t}$ or $\widetilde{R}'_{n,t}$ and, in particular, the fact that (47) gives no rates. Moreover, there is nothing like Theorem 10 here because we have no moment estimates like (31) or (32) for the variables $\widetilde{V}(s)_n$ of (48) (such estimates seem to be out of reach because of the denominators).

**Remark 12.** One could also think of estimators similar to (35), for example,

$$\widetilde{R}''_{n,t} = \sup(r \in \{1,\ldots,d\} : \overline{L}(r)_t^n \geq \rho_n t^{-1/(r-1)}(\overline{L}(r-1)_t^n)^{r/(r-1)}). \qquad (49)$$

However, the previous analysis does *not* carry over to $\widetilde{R}''_{n,t}$, again because we do not know whether the variables $\widetilde{V}(s)_t^n$ converge stably in law on the set $\{R_t < s-1\}$ (once more due to the presence of the denominators).

**Remark 13.** Here, again, the problem of choosing the threshold $\rho_n$ is crucial in practice, despite the fact that a relative threshold is – at least – insensitive to the scale. This is somehow illustrated by the numerical experiments conducted in Section 7. For real data, only the experience of the statistician, at this point, can help for the choice of $\rho_n$. We hope, in a future work, to be able to derive some tests, again based on the consideration of the variables $\overline{L}(r)_t^n$, but which do not necessitate an arbitrary threshold. The next section is something of a first, and somewhat incomplete, attempt in this direction.

## 6. A test based on confidence intervals

Finally, we can take full advantage of the fourth comment in the Introduction. Namely, instead of trying to directly evaluate $R_t$, we can try to evaluate the variables $\overline{L}(s)_t$ for all $s = 1, \ldots, d$.

In view of Theorem 6, this is quite simple. In restriction to the set $\{R_t \geq r\}$, the variables $\sqrt{n}(\overline{L}(r)_t^n - \overline{L}(r)_t)$ are asymptotically mixed normal, with a (conditional) variance $TZ(r,r)_t$ which, in turn, can be estimated by $TZ(r,r)_t^n$ because of Theorem 5. And on the set $\{R_t < r\}$, the variables $\sqrt{n}(\overline{L}(r)_t^n - \overline{L}(r)_t) = \sqrt{n}\overline{L}(r)_t^n$ go to 0 in law.

This allows one to derive (asymptotic) confidence intervals for $\overline{L}(s)_t$. More precisely, we get

$$\lim_n \mathbf{P}(|\sqrt{nTZ(r,r)_t^n}(\overline{L}(r)_t^n - \overline{L}(r)_t)| \geq \gamma \mid R_t = r') = \mathbf{P}(|G| \geq \gamma) \tag{50}$$

for any $\gamma > 0$, where $G$ is an $\mathcal{N}(0,1)$ random variable, provided $\mathbf{P}(R_t = r') > 0$ and $r' \geq r$. This is quite satisfactory because $Z(r,r)_t^n$ is observable. On the other hand, as soon as $\mathbf{P}(R_t = r') > 0$, $r' < r$ and $\gamma > 0$, we get

$$\lim_n \mathbf{P}(|\sqrt{n}\overline{L}(r)_t^n| \geq \gamma \mid R_t = r') = 0. \tag{51}$$

This is less satisfactory because the confidence intervals based on this are not sharp. It is also difficult to obtain non-trivial limit theorems for the sequence $\overline{L}(r)_t^n$, suitably normalized, when $R_t < r$. It seems linked to the speed with which the eigenspaces associated with the positive eigenvalues of $c_s$ rotate in $\mathbf{R}^d$.

***An example.*** Let us consider, for instance, the problem of 'testing' whether the scale-invariant variable $S_t := t^{r-1}\overline{L}(r)_t/(\overline{L}(1)_t)^r$ exceeds some prescribed level $\varepsilon > 0$ (with, say, $0/0 = 0$). This variable is naturally estimated by $S_{n,t} = t^{r-1}\overline{L}(r)_t^n/(\overline{L}(1)_t^n)^r$.

Although, once again, this is not a testing problem in the usual sense, one can proceed as if $S_t$ were a (deterministic) parameter and the null hypothesis is $S_t \geq \varepsilon$. Critical regions on which we reject this hypothesis are naturally of the form

$$C_{n,t}(\eta) = \{S_{n,t} < \eta\}. \tag{52}$$

The 'level' of this test is

$$\alpha_{n,t}^\eta = \sup_{x \geq \varepsilon} \mathbf{P}(C_{n,t}(\eta) \mid S_t \geq x) \tag{53}$$

and its 'power function' for $x \in (0, \varepsilon)$ is

$$\beta_{n,t}^\eta(x) = \mathbf{P}(C_{n,t}(\eta) \mid S_t \leq x) \tag{54}$$

(it would perhaps be more suitable to use $\mathbf{P}(C_n(\eta) \mid S = x)$ as the power function, but the better cannot be evaluated properly below). We also need the variable

$$Z_{n,t} = Tt^{2(r-1)}\frac{(\overline{L}(r)_t^n)^2}{\overline{L}(1)_n^{2r}}\left(\frac{Z(r,r)_t^n}{(\overline{L}(r)_t^n)^2} - \frac{2rZ(1,r)_t^n}{\overline{L}(r)_t^n\overline{L}(1)_t^n} + \frac{r^2Z(1,1)_t^n}{(\overline{L}(1)_t^n)^2}\right), \tag{55}$$

which looks complicated, but is actually computable at stage $n$ from our observations. We then get the following result.

**Theorem 14.** *Assume* (H1) *or* (H2). *Let $\alpha \in (0,1)$ and take $\gamma \in \mathbf{R}$ to be such that $\mathbf{P}(G > \gamma) = \alpha$, where $G$ is an $\mathcal{N}(0,1)$ variable.*

(i) *If $\mathbf{P}(S_t \geq \varepsilon) > 0$, the 'tests' with critical regions $C_{n,t}(\eta_{n,t})$ have an asymptotical level less than or equal to $\alpha$ (i.e., $\limsup_n \alpha_{n,t}^{\eta_{n,t}} \leq \alpha$) if we take*

$$\eta_{n,t} = \varepsilon - \gamma \frac{\sqrt{|Z_{n,t}|}}{\sqrt{n}}. \tag{56}$$

(ii) *If $\mathbf{P}(\overline{L}(1)_t > 0) = 1$, the power function $\beta_{n,t}^{\eta_{n,t}}$ of the above test, with $\eta_{n,t}$ given by (56), satisfies $\beta_{n,t}^{\eta_{n,t}}(x) \to 1$ for any $x \in (0,\varepsilon)$.*

The assumption $\mathbf{P}(\overline{L}(1)_t > 0) = 1$ in (ii) is very mild: it rules out the case where the function $s \mapsto c_s(\omega)$ vanishes on $[0,t]$ on a subset of $\Omega$ with positive probability. If it fails, then the variable $S_t$ is not well defined on this set anyway and the problem is essentially meaningless.

**Proof of Theorem 14.** The result is based on the following consequences of Theorems 5 and 6. We fix $t \in (0,T]$ and introduce the variables

$$\widehat{V}_n := \sqrt{n}(S_{n,t} - S_t) = t^{r-1} \sqrt{n} \left( \frac{\overline{L}(r)_t^n}{(\overline{L}(1)_t^n)^r} - \frac{\overline{L}(r)_t}{(\overline{L}(1)_t)^r} \right). \tag{57}$$

Theorem 6 yields that, in restriction to the set $A = \{\overline{L}(1)_t > 0\}$, the variables $\widehat{V}_n$ converge stably in law to the variable

$$\widehat{V} = t^{r-1} \frac{\overline{L}(r)_t}{(\overline{L}(1)_t)^r} \left( \frac{V(r)_t}{\overline{L}(r)_t} - r \frac{V(1)_t}{\overline{L}(1)_t} \right).$$

Conditionally on the $\sigma$-field $\mathcal{F}$, and again in restriction to $A$, the variable $\widehat{V}$ is centered normal with variance

$$Z = T t^{2(r-1)} \frac{\overline{L}(r)^2}{\overline{L}(1)^{2r}} \left( \frac{Z(r,r)_t}{(\overline{L}(r)_t)^2} - \frac{2r Z(1,r)_t}{\overline{L}(r)_t \overline{L}(1)_t} + \frac{r^2 Z(1,1)_t}{(\overline{L}(1)_t)^2} \right).$$

Finally, the above variable $Z$ is the limit in probability of the sequence $Z_n$ defined in (55), by virtue of Theorem 5. To summarize, we deduce that the variables $T_n = \widehat{V}_n / \sqrt{|Z_{n,t}|}$ converge stably in law, again in restriction to $A$, to an $\mathcal{N}(0,1)$ variable, say $G$, which is independent of $\mathcal{F}$. In particular, for all $y \in \mathbf{R}$, we have

$$B \in \mathcal{F}, B \subset A \implies \mathbf{P}(\{T_n < y\} \cap B) \to \mathbf{P}(B)\mathbf{P}(G < y). \tag{58}$$

It remains to observe that $S_{n,t} = S_t + T_n\sqrt{Z_{n,t}/n}$. With the choice (56) for $\eta_{n,t}$, we then have $T_n < -\gamma$ on $C_{n,t}(\eta_{n,t})$ as soon as $S_t \geq \varepsilon$, and (i) follows from (58) applied to $B = \{S_t \geq \varepsilon\}$, which is included into $A$. Finally, the assumption in (ii) is that $\mathbf{P}(A) = 1$. Let $y, z > 0$ and $x < \varepsilon$, and observe that if $T_n < y$, $S \leq x < \varepsilon$ and $\sqrt{|Z_{n,t}|} \leq z\sqrt{n}$, we have $S_n < x + yz$ and so we are in $C_{n,t}(\eta_{n,t})$ as soon as $yz < \varepsilon - x$. Then (58) with $B = \{S \leq x\}$ and the fact that $Z_{n,t}$ converges in probability to $Z$ yield, for any $y, z > 0$ with $yz < \varepsilon - x$,

$$\mathbf{P}(C_{n,t}(\eta_{n,t}) \cap B) \geq \mathbf{P}(\{T_n < y\} \cap B) - \mathbf{P}(Z_{n,t} > nz^2) \to \mathbf{P}(B)\mathbf{P}(G < y).$$

Since $y$ is arbitrarily large, we conclude that $\mathbf{P}(C_{n,t}(\eta_{n,t}) \cap B) \to \mathbf{P}(B)$ and (ii) follows. $\square$

## 7. Numerical experiments

In this section, we present numerical results for three different families of models. The first two concern financial applications, namely, the calibration of baskets of standard stock prices, or energy indices, with stochastic volatilities. Our last example is not motivated by finance; however, it presents a kind of degeneracy which illustrates an (unsurprising) limitation of our estimation procedure.

All numerical results below concern the test based on a relative threshold described in Section 5.4.

### 7.1. Models with stochastic volatilities

In finance, the calibration of models is a difficult issue. One has to handle missing data in statistical analysis; the frequency of price observations is often too weak to allow one to estimate quadratic variations, and thus volatilities, with good accuracies, and if it is high, then microstructure noises tend to blur the picture; moreover, because of market instabilities, any particular model with fixed parameters or coefficients can pretend to describe market prices over only short periods of time. Consequently, the practitioners are used to calibrating *implicit* parameters of their stock price models by solving PDE inverse problems (see, e.g., Achdou and Pironneau [1] and the references therein), minimizing entropies (see, e.g., Avellaneda *et al.* [2]), etc. Such procedures use instantaneous market information on the stocks under consideration, particularly derivative prices, rather than historical data. In all these approaches, as soon as one deals with a portfolio with several assets, the Brownian dimension is a parameter of prime importance. We have thus studied the performances of our estimation procedure within the commonly used Black, Scholes and Samuelson framework with stochastic volatilities.

Consider

$$\left.\begin{aligned}X_t^1 &= 1 + r_1 \int_0^t X_s^1 \, ds + \sigma_1 \int_0^t X_s^1 \, dB_s^1, \\ X_t^2 &= 1 + r_2 \int_0^t X_s^2 \, ds + \sigma_2 \int_0^t X_s^2 (\rho \, dB_s^1 + \sqrt{1-\rho^2} \, dB_s^2),\end{aligned}\right\}$$

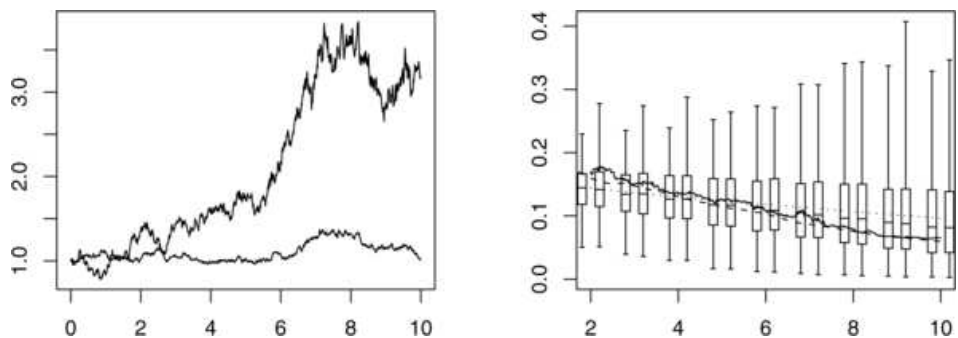

**Figure 1.** Test case 7.1: $\rho = 0$.

with $\sigma_1 = 0.1$, $\sigma_2 = 0.2$, $r_1 = 0.05$ and $r_2 = 0.15$.

To simulate paths of $(X_t^1, X_t^2)$, we have used the Euler scheme with stepsize $10^{-4}$. The final time is $T = 10$. The observations are at times $k \cdot 10^{-2}, 1 \leq k \leq 1000$. In view of (45), we consider the estimators

$$\overline{\xi}_t^n(1) = t\frac{\overline{L}_t^n(2)}{\overline{L}_t^n(1)^2}, \qquad \overline{\xi}_t(1) = t\frac{\overline{L}_t(2)}{\overline{L}_t(1)^2}.$$

Figures 1–4 are organized as follows: the left picture, displays a particular sample path of the pair $(X_t^1, X_t^2)$, and in the right picture, we have plotted the paths of $\overline{\xi}_t^n(1)$ (solid line) and $\overline{\xi}_t(1)$ (dashed line) corresponding to the path of the left display. Moreover, at each integer time $t = 2, 3, \ldots, 10$, the right picture also displays two boxes and whiskers: the box and whiskers on the right plots the empirical quartiles and extends upward and downward to the extremal values of 500 independent samples of the random variables $\overline{\xi}_t^n(1)$; the box and whiskers on the left provides similar information on $\overline{\xi}_t(1)$. Moreover, the left-hand paths in all Figures 1–4 correspond to the same simulated path of the Brownian motion $(B^1, B^2)$.

In this example, we see that the paths of $(X^1, X^2)$ for different values of $\rho$ (and corresponding to the same Brownian motion $(B^1, B^2)$) are difficult to distinguish, whereas the values taken by $\overline{\xi}_t^n(1)$ clearly allow one to distinguish the strongly correlated and weakly correlated cases.

Figures 5 and 6 display the same box and whiskers pictures as previously, but with $\rho = 0.00$ and $\rho = 0.99$, and for two sampling frequencies $T/n$: for each integer time $t$, the left box and whiskers are the same on the left and right displays (they are both for $\overline{\xi}_t(1)$ – beware the change of scale between the two displays), but, unsurprisingly, the spread of $\overline{\xi}_t^n(1)$ is bigger at low frequency (right display). However, even at the lowest frequency (with only 100 observations), it allows correct estimation of the real Brownian dimension.

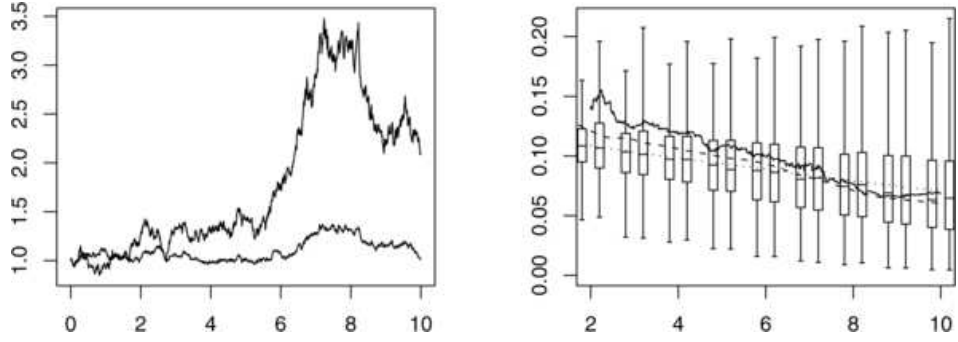

**Figure 2.** Test case 7.1: $\rho = 0.5$.

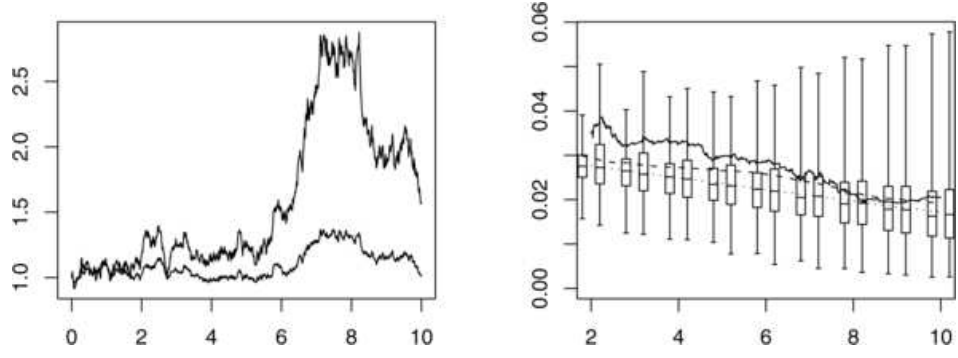

**Figure 3.** Test case 7.1: $\rho = 0.9$.

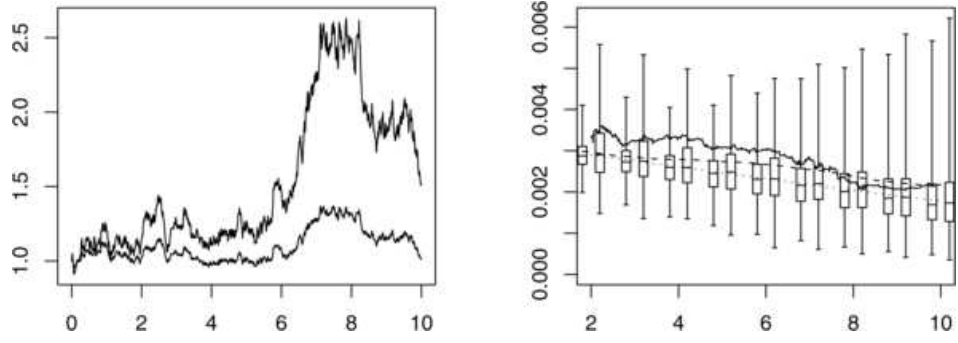

**Figure 4.** Test case 7.1: $\rho = 0.99$.

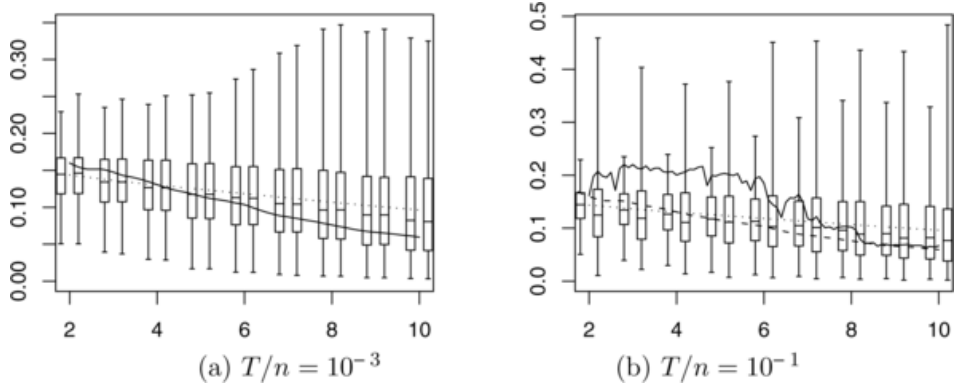

**Figure 5.** Test case 7.1: $\overline{\xi}_t^n(1)$ in terms of $T/n$ ($\rho = 0.00$).

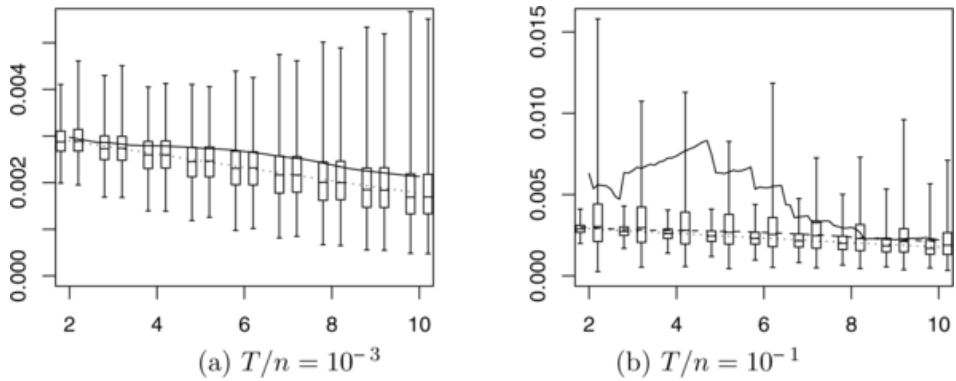

**Figure 6.** Test case 7.1: $\overline{\xi}_t^n(1)$ in terms of $T/n$ ($\rho = 0.99$).

## 7.2. A simplified model for energy indices

We now present a toy model for oil prices. In his Ph.D. thesis (within a collaboration between INRIA and Gaz de France), O. Bardou [3] has studied modelling and simulation questions related to energy contract pricing problems. One question was to identify the coefficients of a stochastic differential system which could satisfyingly describe the dynamics of about ten energy indices.

Here, for the sake of simplicity, we consider a three-dimensional system whose coefficients resemble those identified by O. Bardou: for $1 \leq i \leq 3$, we set

$$dX_t^i = [\alpha_i(X_t^i - K_i)^+ + \beta_i] dB_t^i + \nu_i(\mu_i - X_t^i) dt,$$

with $X_0^1 = 0.29$, $X_0^2 = 0.89$, $X_0^3 = 0.62$.

Fixed this way, the diffusion term does not satisfy Hypothesis (H2). We thus slightly modify the equation and consider

$$dX_t^i = [\alpha_i \phi(X_t^i - K_i) + \beta_i] \, dB_t^i + \nu_i(\mu_i - X_t^i) \, dt,$$

where $\phi(x) = 0$ if $x \leq 0$, $\phi(x) = 2.5x^2$ if $0 < x < 0.2$ and $\phi(x) = x - 0.1$ if $x \geq 0.2$. In this very simplified model, the components of $X$ are independent; in the real situation where one observes energy indices, one should take correlated Brownian motions $B^i$, as in Section 7.1.

We set $\nu_i = \mu_i = \alpha_i = 1$. The drift term then stabilizes the process around the value 1. If $\beta_i = 0$, the process $(X_t^i)$ diffuses only when $X_t^i$ is above the threshold $K_i$.

As above, we have approximated a path of the system by simulating the Euler scheme with stepsize $10^{-4}$ between times 0 and 10. The observations are at times $k \cdot 10^{-2}, 1 \leq$

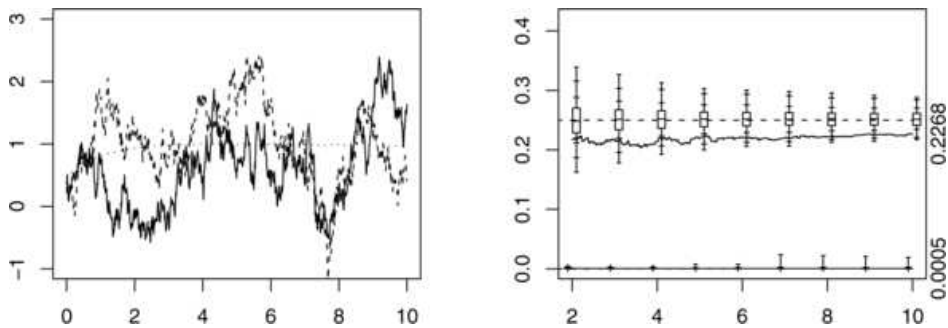

**Figure 7.** Test case 7.2: $\beta_1 = \beta_2 = 1$, $\beta_3 = 0$, $K_1 = K_2 = 3$, $K_3 = 0.9$.

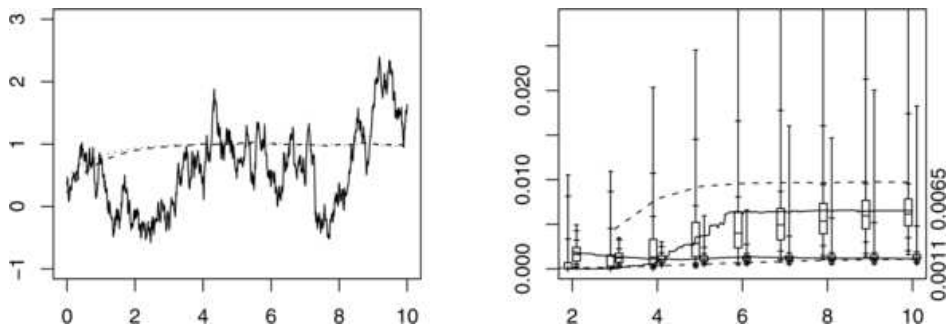

**Figure 8.** Test case 7.2: $\beta_1 = 1$, $\beta_2 = \beta_3 = 0$, $K_1 = 3$, $K_2 = K_3 = 0.9$.

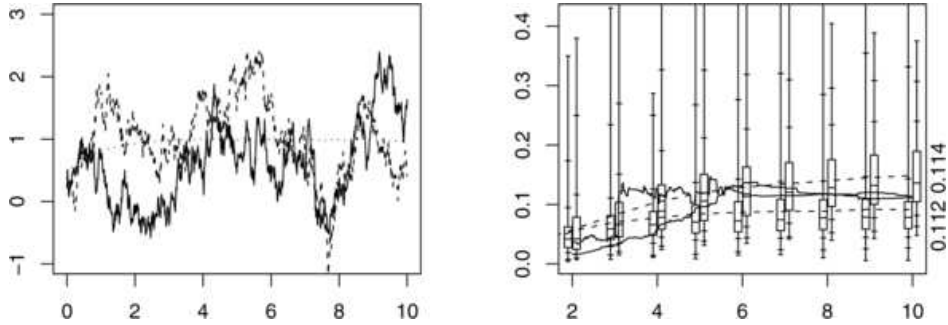

**Figure 9.** Test case 7.2: $\beta_1 = 1$, $\beta_2 = \beta_3 = 0$, $K_1 = 3$, $K_2 = K_3 = 0.6$.

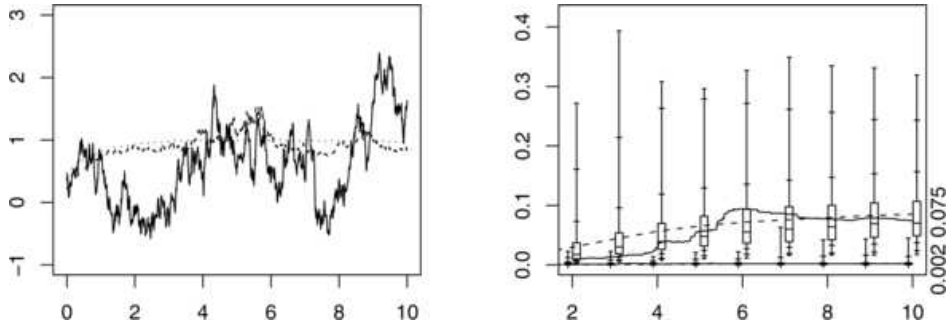

**Figure 10.** Test case 7.2: $\beta_1 = 1$, $\beta_2 = \beta_3 = 0$, $K_1 = 3$, $K_2 = 0.6$, $K_3 = 0.9$.

$k \leq 1000$. In view of (45), we consider the estimators

$$\overline{\xi}_t^n(1) = t\frac{\overline{L}_t^n(2)}{\overline{L}_t^n(1)^2}, \qquad \overline{\xi}_t^n(2) = \sqrt{t}\frac{\overline{L}_t^n(3)}{\overline{L}_t^n(2)^{3/2}}$$

and

$$\overline{\xi}_t(1) = t\frac{\overline{L}_t(2)}{\overline{L}_t(1)^2}, \qquad \overline{\xi}_t(2) = \sqrt{t}\frac{\overline{L}_t(3)}{\overline{L}_t(2)^{3/2}}.$$

In Figures 7–10, the left boxes show a particular path of $(X^1, X^2, X^3)$. The right boxes show the corresponding paths of $\overline{\xi}_t^n(j)$ (solid line) and $\overline{\xi}_t(j)$ (dashed line) for $j = 1$ and $j = 2$, respectively the top and the bottom curves. Values on the right-hand side vertical axes denote $\overline{\xi}_T^n(1)$ and $\overline{\xi}_T^n(2)$. Moreover, on the right, we have box and whiskers for the empirical quartiles of $\overline{\xi}_t^n(1)$ (top) and $\overline{\xi}_t^n(2)$ (bottom), computed from 500 independent paths and for all integer times $t = 2, 3, \ldots, 10$. The whiskers extend to the extremal values of the samples, with the other ticks denoting the 1%, 10%, 90% and 99% quantiles.

In Figure 7, the two first components diffuse from time 0 to time 10 since $K_1$ and $K_2$ are large. The third component diffuses only a little since it is attracted to 1 and $\phi(1-K_3) = \phi(0.1)$ is small. Given the threshold $\rho_n = 0.01$, the explicative Brownian dimension $\widetilde{R}_{n,t}$ is 2 since $\overline{\xi}^n_t(1)$ takes values around 0.2, whereas $\overline{\xi}^n_t(2)$ takes values around $2 \cdot 10^{-3}$.

In Figure 8, the two last components have a small diffusion term. As both $\overline{\xi}^n_t(1)$ and $\overline{\xi}^n_t(2)$ take values less than 0.02, according to the same threshold $\rho_n = 0.01$ as above, the explicative Brownian dimension $\widetilde{R}_{n,t}$ is 1.

In Figure 9, we have $K_1$ large and $\phi(1-K_2) = \phi(1-K_3) = \phi(0.4) = 0.3$. Therefore, none of the diffusion terms can be neglected, but the first component 'oscillates' more than the two others. It is a case where observed paths, for which both $\overline{\xi}^n_t(1)$ and $\overline{\xi}^n_t(2)$ take values around 0.1, may make it difficult to decide whether the explicative Brownian dimension should be chosen as 1 or 2.

Finally, in Figure 10, we keep the first two components as in Figure 9, but change the third into an almost constant process. Of course, as above, we have difficulties in deciding whether the Brownian dimension is 1 or 2. However, it is clear that it cannot be 3 since $\overline{\xi}^n_t(2)$ fluctuates around $5 \cdot 10^{-3}$.

## 7.3. Sensitivity to a drift term close to be a martingale

We now consider a model with a strongly oscillating drift term. These oscillations significantly decrease the efficiency of our estimator. In certain circumstances, the explicative Brownian dimension overestimates the real dimension.

The system under consideration is

$$\left. \begin{aligned} X^1_t &= \int_0^t \eta \cos(\theta X^2_s) \, ds, \\ X^2_t &= B_t, \end{aligned} \right\},$$

where $B$ is a one-dimensional standard Brownian motion and $\eta$, $\theta$ are positive real numbers.

As above, we have approximated a path of $X$ by simulating the Euler scheme with stepsize $10^{-4}$ between times 0 and 10. The observations are at times $k \cdot 10^{-2}, 1 \leq k \leq 1000$. In view of (45), we consider the estimator $\overline{\xi}^n_t(1) = t\frac{\overline{L}^n_t(2)}{\overline{L}^n_t(1)^2}$. Observe that $\overline{L}_t(2) = 0$, so that $\overline{\xi}^n_t(1)$ should be close to 0.

The boxes and whiskers denote the empirical quartiles of $\overline{\xi}^n_t(1)$ and extend to the extremal values, at times $t = 2, 3, \ldots, 10$ and for 500 estimations of $\overline{\xi}^n_t(1)$.

Figure 11 shows a case with such a highly oscillating coefficient that the first component is close to 0. The reason is clear since

$$\begin{aligned} X^1_t &= \frac{2\eta}{\theta^2} - \frac{2\eta}{\theta^2} \cos(\theta B_t) \\ &\quad - \frac{2\eta}{\theta} \int_0^t \sin(\theta B_s) \, dB_s \end{aligned} \tag{59}$$

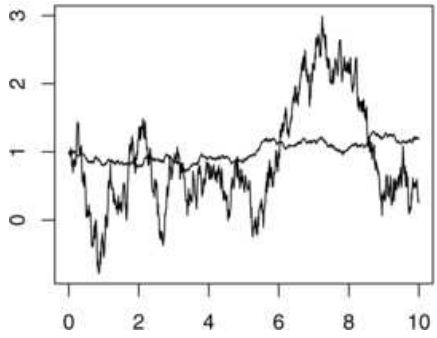 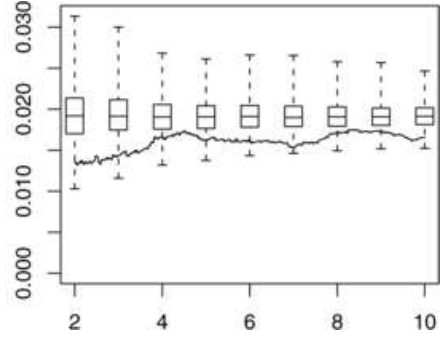

**Figure 11.** Test case 7.3: $\eta = 10$, $\theta = 100$.

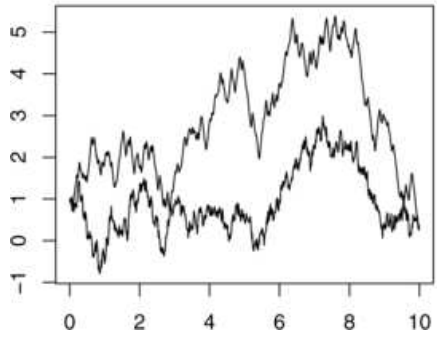 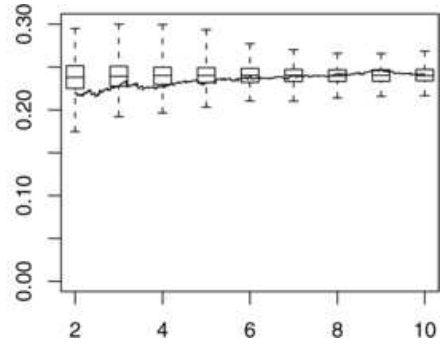

**Figure 12.** Test case 7.3: $\eta = 10$, $\theta = 10$.

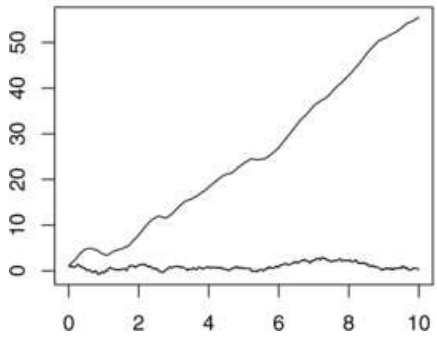 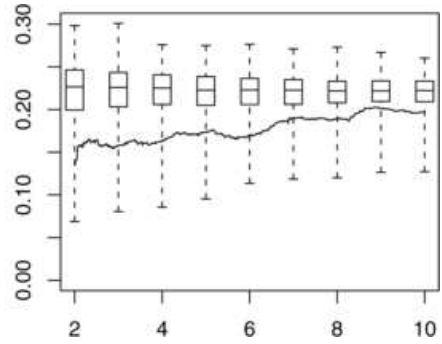

**Figure 13.** Test case 7.3: $\eta = 10$, $\theta = 1$.

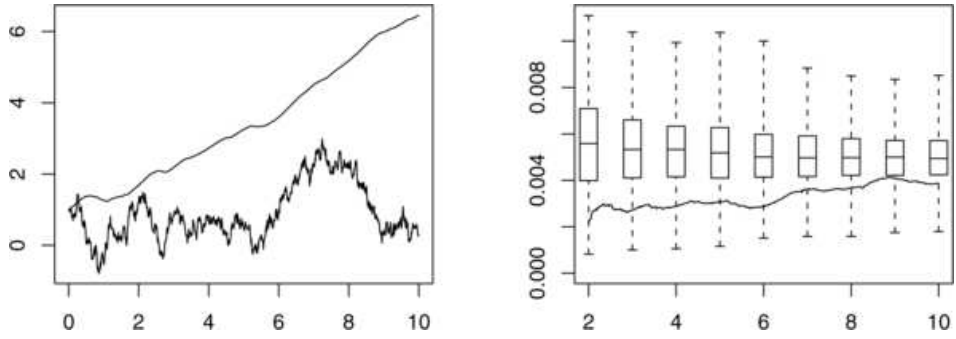

**Figure 14.** Test case 7.3: $\eta = 1$, $\theta = 1$.

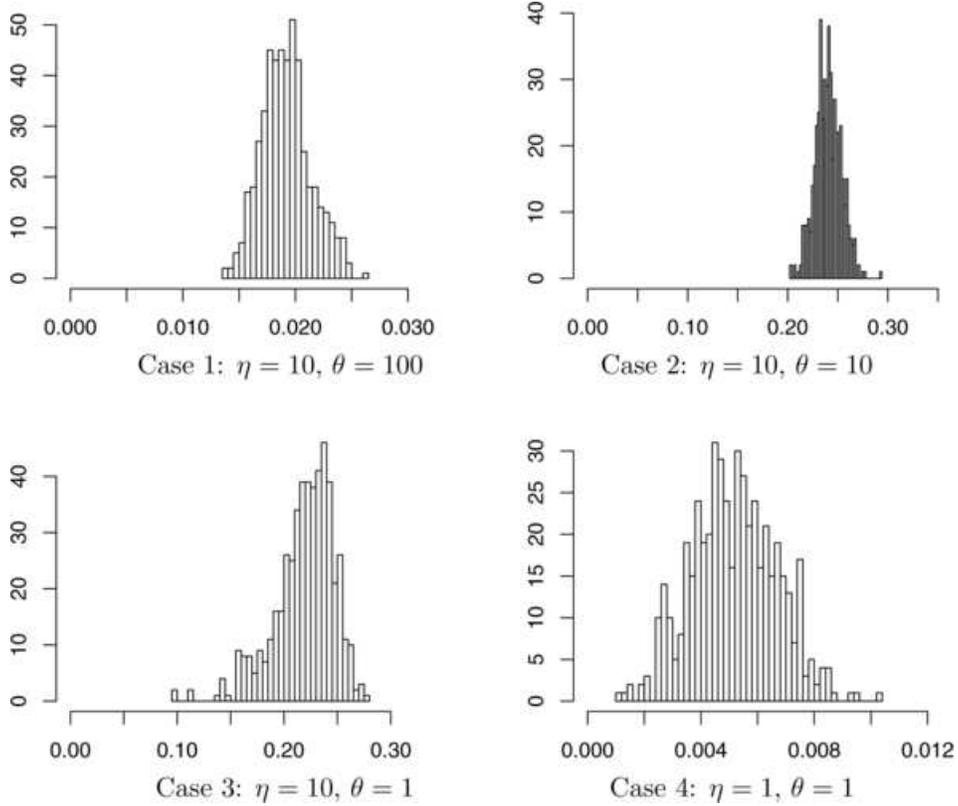

**Figure 15.** Test case 7.3: Histograms of $\overline{\xi}_t^n(1)$ at time $t = 5$.

and, here, $\eta = 10$ and $\theta = 100$. As our estimator takes values around 0.02, choosing $\rho_n = 0.02$ leads one to choose $\widetilde{R}_{n,t} = 1$ as the explicative Brownian dimension.

In Figure 12, we fix $\eta = 10$ and $\theta = 10$. In view of (59), $(X_t^1)$ is close to being a stochastic integral whose covariation with $B$ on the time interval $[0, 10]$ is small. As our estimator now takes values around 0.2, we are led to choose 2 as the explicative Brownian dimension and thus to overestimate the real Brownian dimension.

In Figure 13 we fix $\eta = 10$ and $\theta = 1$. As our estimator takes values around 0.2, we again overestimate the real Brownian dimension.

In Figure 14, as $\eta = 1$ and $\theta = 1$, the first component oscillates 'reasonably'. The estimator takes values less than 0.01 and we are led to correctly choose 1 as the explicative Brownian dimension.

Figure 15 shows histograms of $\overline{\xi}_t^n(1)$ in the various preceding situations. Finally, Figure 16 shows the influence of the sampling frequency in the case exhibited in Figure 14: we see that in this case, the Brownian dimension remains correctly estimated as 1 when the step size remains smaller than $10^{-1}$.

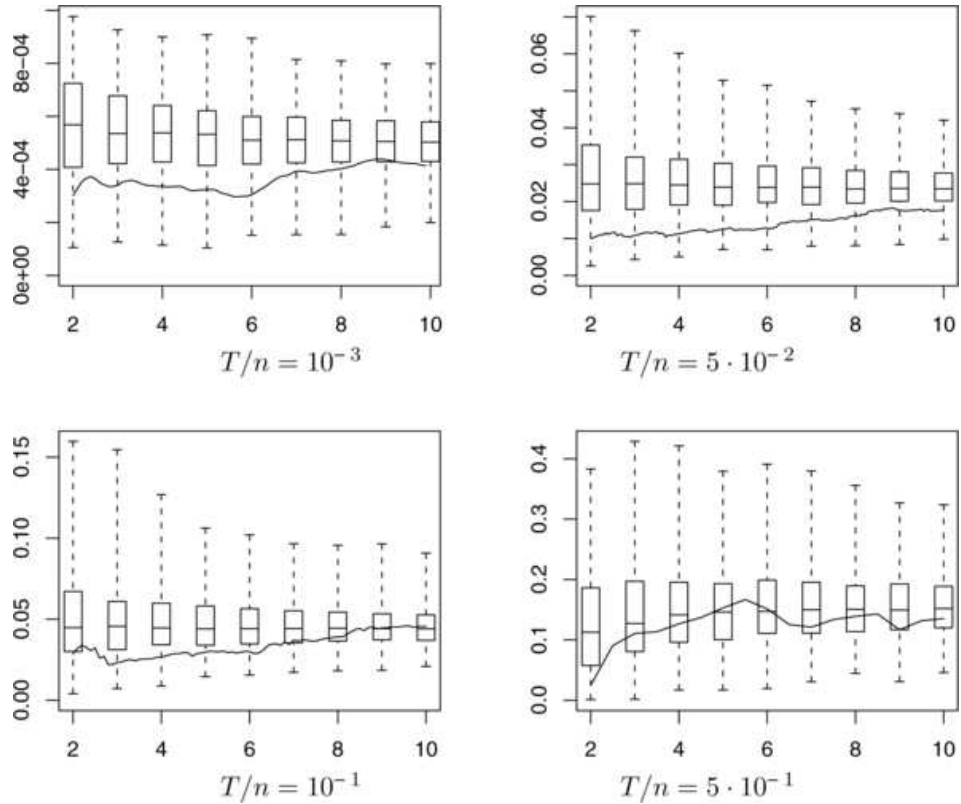

**Figure 16.** Test case 7.3: $\overline{\xi}_t^n(1)$ in terms of $T/n$ ($\eta = 1$, $\theta = 1$).